\documentclass[11pt,leqno]{article}
\usepackage{graphicx, amsfonts, amsthm, amsxtra, verbatim, multicol}
\usepackage[mathscr]{euscript}
\begin{document}

\begin{center}
{\LARGE\bf
Calculation of mixed Hodge structures, Gauss-Manin connections and 
Picard-Fuchs equations
}
\footnote{ Math. classification: 14C30,
32S35 
} \\
\vspace{.25in} {\large {\sc Hossein Movasati}} \\
Technische Universit\"at Darmstadt \\
Fachbereich Mathematik \\
Schlo\ss gartenstr. 7, 
64289 Darmstadt, Germany \\ 
Email: {\tt movasati@mathematik.tu-darmstadt.de} 
\end{center}

\begin{abstract}
In this article we introduce algorithms which compute 
iterations of Gauss-Manin connections, Picard-Fuchs equations of Abelian integrals and
mixed Hodge structure of affine varieties of dimension $n$ 
in terms of differential forms.
In the case $n=1$ such computations have many applications in differential
equations and counting their limit cycles. For $n>3$, these computations 
give us an explicit definition of Hodge cycles. 
\end{abstract}

\newtheorem{theo}{Theorem}
\newtheorem{exam}{Example}
\newtheorem{coro}{Corollary}
\newtheorem{defi}{Definition}
\newtheorem{prob}{Problem}
\newtheorem{lemm}{Lemma}
\newtheorem{prop}{Proposition}
\newtheorem{rem}{Remark}
\newtheorem{conj}{Conjecture}
\def\hol{{\mathrm Hol}}        
\def\sing{{\mathrm Sing}}            
\def\cha{{\mathrm char}}             
\def\Gal{{\mathrm Gal}}              
\def\jacob{{\mathrm jacob}}          
\newcommand\Pn[1]{\mathbb{P}^{#1}}   
\def\Z{\mathbb{Z}}                   
\def\Q{\mathbb{Q}}                   
\def\C{\mathbb{C}}                   
\def\R{\mathbb{R}}                   
\def\N{\mathbb{N}}                   
\def\A{\mathbb{A}}                   
\def\uhp{{\mathbb H}}                
\newcommand\ep[1]{e^{\frac{2\pi i}{#1}}}
\newcommand\HH[2]{H^{#2}(#1)}        
\def\Mat{{\mathrm Mat}}              
\newcommand{\mat}[4]{
     \begin{pmatrix}
            #1 & #2 \\
            #3 & #4
       \end{pmatrix}
    }                                
\newcommand{\matt}[2]{
     \begin{pmatrix}                 
            #1   \\
            #2
       \end{pmatrix}
    }
\def\ker{{\mathrm ker}}              
\def\cl{{\mathrm cl}}                
\def\dR{{\mathrm dR}}                
\def\Hb{{\cal H}}                    
\def\GL{{\mathrm GL}}                
\def\pedo{{\cal P}}                  
\def\PP{\tilde{\cal P}}              
\def\cm {{\cal C}}                   
\def\K{{\mathbb K}}                  
\def\F{{\cal F}}                     
\def\M{{\cal M}}
\def\RR{{\cal R}}
\newcommand\Hi[1]{\mathbb{P}^{#1}_\infty}
\def\pt{\mathbb{C}[t]}               
\def\W{{\cal W}}                     
\def\Af{{\cal A}}                    
\def\gr{{\mathrm Gr}}                
\def\Im{{\mathrm Im}}                
\newcommand\SL[2]{{\mathrm SL}(#1, #2)}    
\newcommand\PSL[2]{{\mathrm PSL}(#1, #2)}  
\def\Res{{\mathrm Res}}              

\def\L{{\cal L}}                     
\def\Aut{{\mathrm Aut}}              
\def\any{R}                          
\newcommand\ovl[1]{\overline{#1}}    

\section{Introduction}
\label{intro}
The theory of Abelian integrals which arises in polynomial 
differential equations of the
type $\dot x=P(x,y),\ \dot y=Q(x,y)$ is one of 
the most fruitful areas which needs  a special attention form algebraic 
geometry. The reader is referred 
to the articles \cite{Il02}, \cite{mov0} and  \cite{ga02} 
for a history and applications 
of such Abelian integrals in differential equations. 
In this article we deal with computational aspects of such integrals. 
All polynomial objects which we use are defined over $\C$.

Let us be given a polynomial $f$ in $n+1$ variables
$x_1,x_2,\ldots, x_{n+1}$, a polynomial differential $n$-form $\omega$ and a 
continuous family of $n$-dimensional oriented cycles 
$\delta_t\subset L_t:=f^{-1}(t)$. The 
protagonist of this article is the integral $\int_{\delta_t}\omega$, called 
Abelian integral. 
Computations related to these integrals becomes easier when we put a certain
kind of tameness condition on $f$ (see \S \ref{adami}). For such a tame 
polynomial we can write $\int_{\delta_t}\omega$ as:
\begin{equation}
\label{26nov}
\sum_{\beta\in I}p_\beta(t)\int_{\delta_t}\eta_\beta
\end{equation}
where $\eta_\beta,\beta\in I$ is a class of differential $n$-forms constructed
from a basis of the Milnor vector space of $f$ and $p_\beta$'s are polynomials 
in $t$ (see \S \ref{brieskorn} for the algorithm which produces $p_\beta$'s).
The Guass-Manin connection $\nabla\omega$ has the following basic
property 
\begin{equation}
\label{gofta}
\frac{\partial}{\partial t}\int_{\delta_t}\omega=\int_{\delta_t}
\nabla\omega
\end{equation}
The above term can be written in the form (\ref{26nov}) with $p_\beta$'s 
rational functions in $t$ with poles in the critical values of $f$ (see \S 
\ref{gaussmanin} for the algorithm which produces $p_\beta$'s). 
The $n$-th cohomology  of a smooth fiber $L_t$ is canonically 
isomorphic to $\Omega_{L_t}^n/\Omega_{L_t}^{n-1}$, where 
$\Omega_{L_t}^i$ is the restriction  of polynomial differential $i$-forms 
to $L_t$, 
and carries two natural filtrations called the weight and the Hodge filtrations
(both together is called the mixed Hodge structure).
These  filtrations are generalizations of classical notions of differential 
forms of the first, second and third type for Riemann surfaces in 
higher dimensional varieties. The reader who is not interested in the case
$n>1$ is invited to follow the article with $n=1$ and with the usual
notions of differential forms of the first, second and third type. 
How to calculate these filtrations by
means of differential forms is the main theorem of \cite{mo} and related
algorithms are explained in \S \ref{dbeta}.  Last but not the least, 
our protagonist satisfies a Picard-Fuchs equation 
$\sum_{i=0}^kp_i(t)\frac{\partial^i}{\partial t^i}=0$, where $p_i$'s are 
polynomials in $t$. The algorithm which produces $p_\beta$'s is explained
in \S \ref{pfeq}. The theory of Abelian integrals can be studied even in the case $n=0$, i.e. $f$ is a polynomial in one variable. Since some open problems,
for instance infinitesimal 
Hilbert Problem (see \cite{Il02}), can be also stated in 
this case, we have included \S \ref{n=0}. 
All the algorithms explained in this article are
implemented in the library {\tt brho.lib} of {\sc Singular} and throughout
the article we mention the related procedures of {\tt brho.lib}. 
Nevertheless, the reader may want to implement the algorithms of this article
in any other software in commutative algebra.

The main theorem of \cite{mo}(see also \S \ref{maintheo}) does not give 
a basis of the Brieskorn module compatible with the mixed Hodge structure
(see Definition \ref{mhs}). In \S \ref{examples} we 
obtain such bases for some examples of $f$ by modifying the one given in
\S \ref{maintheo} (we do not have a general method for every $f$). 
Applications of our computations in differential equations and particularly 
in direction of the  article \cite{ga02} is a matter of future work.

\section{Tame polynomials and Brieskorn modules}
\label{adami}
We start with a definition.
\begin{defi}\rm
\index{tame polynomial}
A polynomial $f\in\C[x]$ is called (weighted) tame if 
there exist natural numbers  $\alpha_1,\alpha_2,\ldots,\alpha_{n+1}\in\N$ 
such that
${\mathrm Sing}(g)=\{0\}$, where $g=f_d$ is the last 
homogeneous piece of $f$ in the graded algebra 
$\C[x],\ {\mathrm deg}(x_i)=\alpha_i$. \index{graded ring}
\end{defi}
The multiplicative group $\C^*$ acts on
$\C^{n+1}$ in the following way:
$$
\lambda^*:(x_1,x_2,\ldots,x_{n+1}) \rightarrow
(\lambda^{\alpha_1}x_1,\lambda^{\alpha_2}x_2,
\ldots,\lambda^{\alpha_{n+1}}x_{n+1}),\ \lambda\in\C^*
$$
The polynomial (resp. the polynomial form) 
$\omega$ in $\C^{n+1}$ is (weighted) homogeneous of degree $d\in\N$ if
$\lambda^*(\omega)=\lambda^d \omega, \ \lambda\in\C^*$.
Fix a homogeneous polynomial $g$ of degree $d$ and with an isolated 
singularity at $0\in \C^{n+1}$. Let $\Af_g$ be the affine space
of all tame polynomials $f=f_0+f_1+\cdots+f_{d-1}+g$. The space $\Af_g$ is
parameterized by the coefficients of $f_i, i=0,1,\ldots,d-1$.
The multiplicative group
$\C^{*}$ acts on $\Af_g$ by
$$
\lambda\bullet f=\frac{f\circ \lambda^*}{\lambda^d}=\lambda^{-d}f_0+
\lambda^{-d+1}f_1+\cdots+\lambda^{-1}f_d+g
$$
The action of $\lambda\in \C^*$
takes $\lambda\bullet f=0$ biholomorphically to $f=0$.
 
Let $f\in\Af_g$. 
We choose a basis
$x^I:=\{x^\beta\mid\beta\in I\}$ of monomials for the Milnor $\C$-vector space
$$
V:=\C[x]/\jacob(g)
$$
In {\sc Singular} one can get $x^I$ using {\tt kbase} command. Using this
command $\deg(x^\beta)$ is not a decreasing sequence. In {\tt brho.lib} 
the procedure {\tt okbase} makes a permutation on the result of 
{\tt kbase} and gives us $x^\beta$'s with $\deg(x^\beta)$ decreasing.
We will fix the order obtained by {\tt okbase} and use it for many 
examples throughout the text. Define
\begin{equation}
\label{21may2004}
w_i:=\frac{\alpha_i}{d},\ 1\leq i\leq n+1,\ \eta:=(\sum_{i=1}^{n+1}(-1)^{i-1}
w_i x_i\widehat {dx_i}),\ L_t:=f^{-1}(t),t\in\C
\end{equation}
$$
 A_\beta:=
\sum_{i=1}^{n+1}(\beta_i+1)w_i,\
\eta_\beta:=x^\beta \eta,\ \omega_\beta=x^\beta dx \ \beta\in I
$$
where $\widehat{dx_i}=dx_1\wedge\cdots
\wedge dx_{i-1}\wedge dx_{i+1}\wedge\cdots\wedge dx_{n+1}$. 
Note that $A_\beta=\frac{\deg(x^{\beta+1})}{d}$. 
It turns out
that $x^I$ is also a basis of $V_f:=\C[x]/\jacob(f)$ and so $f$ and $g$ have the same Milnor 
numbers (see the conclusion after Lemma 4 of \cite{mo}). 
We denote  it by  $\mu$ . 
We denote by $P$ 
the set of critical points of $f$ and by $C:=f(P)$ the set of critical values 
of $f$. We will use also $P$ for a polynomial in $\C[x]$. This will not
make any confusion!.

Let $\Omega^i, i=1,2,\ldots, n+1$ (resp. $\Omega^i_j,\ j\in\N\cup\{0\}$) 
be the set of polynomial differential $i$-forms (resp. homogeneous
polynomial differential $i$-forms) in $\C^{n+1}$. The Milnor vector space of $g$ 
can be rewritten in the form $V:=
\frac{\Omega^{n+1}}{df\wedge \Omega^n}$. 
The Brieskorn modules
$$
H'=H'_f:=\frac
{\Omega^n}{df\wedge \Omega^{n-1}+d\Omega^{n-1}},\ H''=H''_f=
\frac{\Omega^{n+1}}{df\wedge d\Omega^{n-1}}
$$
of $f$ are $\C[t]$-modules
in a natural way: $t.[\omega]=[f\omega],\ [\omega]\in H' $ resp. $\in H''$.
They are defined in the case $n>0$. 
The case $n=0$ is treated separately in \S \ref{n=0}. 
\begin{defi}\rm
Let $H$ be one of $H'$ or $H''$. 
If $H=H''$ then by restriction of  $\omega$ on $L_c,\ c\in\C\backslash C$ 
we mean the residue of 
$\frac{\omega}{f-c}$ in $L_c$ and by $\int_{\delta}\omega,\ \delta\in 
H_n(L_c,\Z)$ we mean $\int_\delta Resi(\frac{\omega}{f-c})$. It is natural to 
define the mixed Hodge structure of $H$ as follows:
$W_mH,\ m\in\Z$ (resp. $F^kH,\ k\in\Z$) 
consists of elements $\omega\in H$ such that the restriction of 
$\omega$ on all $L_c,\ c\in\C\backslash C$ belongs to 
$W_mH^n(L_c,\C)$ (resp. $F^kH^n(L_c,\C)$). 
\end{defi}\rm
In the case of a tame polynomial we have
$$
\{0\}=W_{n-1}\subset W_n \subset W_{n+1}=H,\ \  \{0\}=F^{n+1}\subset F^{n}\subset
\cdots\subset F^1\subset F^0=H
$$
Each piece of the mixed Hodge structure of $H$ is
a $\C[t]$-module. In the same way we define the mixed Hodge structure
of the localization of $H$ over multiplicative subgroups of $\C[t]$.
In the case $n=1$ our definition can be simplified as follows: 
We have the filtrations 
$\{0\}=W_0\subset W_1\subset W_2=H$ and
$0=F^2\subset F^1\subset F^0=H$, where 
$$
W_1=\{\omega \in H\mid \omega \hbox{ restricted to a regular fiber 
has not residue at infinity }  \}
$$
$$
F^1=\{
\omega\in H\mid \omega 
\hbox{ restricted to a regular fiber has poles of maximum order 1 at infinity}
\}
$$
In particular we get
$$
W_1\cap F^1=\{
\omega\in H\mid \omega 
\hbox{ restricted to a regular comactified  fiber is of the first kind}
\}
$$
For the notion of compactification of $\C^2$ and infinity see \cite{ga02} and
\cite{mov}. 

The projection of $F^\bullet$ in $\gr_m^W H:=W_{m}/W_{m-1}$ gives us
the filtration $\bar F^\bullet$ in $\gr_m^W H$ and we define
$\gr^k_F\gr^W_m H=\bar F^{k}/\bar F^{k+1}$.
\begin{defi}
\label{mhs}
Suppose that $H$ is freely generated module. 
The set $B=\cup_{m,k\in \Z} B_m^k\subset H$ is a basis of 
$H$ compatible with the mixed Hodge
structure  if $B_m^k$ form a basis of $\gr^k_F\gr^W_m H$.
\end{defi} 
\section{Quasi-homogeneous singularities}
\label{homogeneous}
Let $f=g$ be a weighted homogeneous polynomial \index{weighted homogeneous
polynomial} with an isolated singularity at origin. We explain the algorithm
which writes every element of  $H''$ of $g$ as a
$\C[t]$-linear combination of $\omega_\beta$'s. Recall that
$$
dg\wedge d(P\widehat{dx_i,dx_j})=(-1)^{i+j +\epsilon_{i,j} }
(\frac{\partial g}{\partial x_j}\frac{\partial P}{\partial x_i}
-\frac{\partial g}{\partial x_i }\frac{\partial P}{\partial x_j})dx
$$
where $\epsilon_{i,j}=0$ if $i<j$ and $=0$ if $i>j$ and 
$\widehat{dx_i,dx_j}$ is $dx$ without $dx_i$ and $dx_j$ (we have not changed
the order of $dx_1,dx_2,\ldots$ in $dx$).
\begin{prop}
For a monomial $P=x^\beta$ we have
\begin{equation}
\label{4.7.03}
\frac{\partial g}{\partial x_i}.Pdx=\frac{d}{d.A_\beta-\alpha_i}
\frac{\partial P}{\partial x_i}gdx+
dg\wedge d(\sum_{j\not=i}\frac{(-1)^{i+j+1+\epsilon_{i,j}}\alpha_j}{dA_\beta-\alpha_i}
x_jP \widehat{dx_i,dx_j})
\end{equation}
\end{prop}
\begin{proof}
The proof is a straightforward calculation.
$$
 \sum_{j\not=i}
\frac{(-1)^{i+j+1+\epsilon_{i,j}}\alpha_j}{dA_\beta-\alpha_i}
dg\wedge d(x_jP \widehat{dx_i,dx_j}) 
$$
$$ 
 \frac{-1}{dA_\beta-\alpha_i}\sum_{j\not=i}(
\alpha_j\frac{\partial g}{\partial x_j}\frac{\partial (x_j P)}{\partial x_i}
-\alpha_j 
\frac{\partial g}{\partial x_i }\frac{\partial (x_j P)}{\partial x_j}) 
$$
$$
=\frac{-1}{dA_\beta-\alpha_i}( (dg- \alpha_ix_i\frac{\partial g}{\partial x_i})
\frac{\partial P}{\partial x_i}-P\frac{\partial g}{\partial x_i}
\sum_{j\not=i} \alpha_j(\beta_j+1))
$$
$$
\frac{-1}{dA_\beta-\alpha_i}(dg\frac{\partial P}{\partial x_i}-\alpha_i\beta_i
P\frac{\partial g}{\partial x_i}-P\frac{\partial g}{\partial x_i}
\sum_{j\not=i} \alpha_j(\beta_j+1) )
$$
\end{proof}
We use the above Proposition to write every $Pdx\in \Omega^{n+1}$ in 
the form 
\begin{equation}
\label{14.11}
Pdx=\sum_{\beta\in I}p_\beta(g)\omega_\beta+
dg\wedge d\xi,\  p_\beta\in\C[t],\ \xi\in\Omega^{n-1},
\deg(p_\beta(g)\omega_\beta,dg\wedge d\xi)\leq \deg(Pdx)
\end{equation}
\begin{itemize}
\item 
Input: The homogeneous polynomial $g$ and $P\in\C[x]$ representing 
$[Pdx]\in H''$

Output: $p_\beta,\beta\in I$ and $\xi$ satisfying (\ref{14.11})

We write 
\begin{equation}
\label{22:51}
Pdx=\sum_{\beta\in I}c_\beta x^\beta.dx+dg\wedge\eta, \ 
\deg(dg\wedge\eta)\leq \deg(Pdx)
\end{equation} 
Then we apply (\ref{4.7.03}) to each monomial component 
$\tilde P\frac{\partial g}{\partial x_i}$  
of $dg\wedge\eta$ and then  we write each 
$\frac{\partial \tilde P}{\partial x_i}dx$ in the form (\ref{22:51}). 
The degree of the components which make $Pdx$ not to be of the form 
(\ref{14.11}) always 
decreases and finally we get the desired form. 
\end{itemize}
This algorithm is implemented  in the procedure {\tt linear1} of 
the library  {\tt brho.lib}. To find  a similar algorithm for
$H'$ we note that  if $\eta\in \Omega^{n}$ is written in the
form
\begin{equation}
\label{29.6}
\eta=\sum_{\beta\in I}p_\beta(g)\eta_\beta+dg\wedge\xi+d\xi_1,\ p_\beta\in
\C[t],\xi,\xi_1\in\Omega^{n-1}  
\end{equation}
where each piece in the right hand side of the above equality has degree less 
than $\deg(\eta)$  then
\begin{equation}
\label{roesler}
d\eta=\sum_{\beta\in I}(p_\beta(g)A_\beta+p'_\beta(g)g)\omega_\beta-
dg\wedge d\xi
\end{equation}
and the 
inverse of the  map $\C[t]\rightarrow\C[t],\ p\rightarrow A_\beta.p+p'.t$
is given by $\sum_{i=0}^ka_it^i\rightarrow \sum_{i=1}^k
\frac{a_i}{A_\beta+i}t^i$. Therefore, we can recover 
$p_\beta$'s using {\tt linear1} for $d\eta$. 
The obtained algorithm is implemented in 
the procedure {\tt linear2}  of the library  {\tt brho.lib}.
Later we will introduce the procedure {\tt linear} (resp. {\tt linearp}) 
which works for arbitrary tame polynomial and is an extended version 
of {\tt linear1} (resp. {\tt linear2}). For this reason, the procedures
{\tt linear1, linear2} are not available to the user.
  
\begin{theo}
\label{stee}
For a weighted homogeneous polynomial $g$, 
the set 
$$
B=\cup_{k=1}^{n}B^k_{n+1}\cup
 \cup_{k=0}^{n} B_{n}^k
$$
with
$$ 
B^{k}_{n+1}=\{\omega_\beta\mid  A_\beta=n-k+1\},
B^{k}_{n}=\{\omega_\beta\mid  n-k<A_\beta<n-k+1\},
$$
is a basis of $H''$ compatible with the mixed Hodge structure.
The same is true for $H'$ replacing $\omega_\beta$ with $\eta_\beta$. 
\end{theo}
This theorem is due to J. Steenbrink \cite{st77} and its generalization for an
arbitrary tame polynomial is given in \S \ref{maintheo}.
\section{A basis of $H'$ and $H''$}
\label{brieskorn}
\begin{prop}
For every tame polynomial $f\in\Af_g$ the forms
$\omega_\beta,\ \beta\in I$ (resp. $\eta_\beta,\ \beta\in I $) form a basis
of the Brieskorn module $H''$ (resp. $H'$) of $f$. More precisely, every 
$\omega\in\Omega^{n+1}$ (resp. $\omega\in \Omega^{n}$) can be written 
\begin{equation}
\label{abelian}
\omega=\sum_{\beta\in I}p_\beta(f)\omega_\beta+
df\wedge d\xi,\  p_\beta\in\C[t],\ \xi\in\Omega^{n-1},\ 
\deg(p_\beta)\leq \frac{\deg(\omega)}{d}-A_\beta
\end{equation}
(resp.
\begin{equation}
\label{abelianp}
\omega=\sum_{\beta\in I}p_\beta(f)\eta_\beta+
df\wedge \xi+d\xi_1,\  p_\beta\in\C[t],\ \xi\in\Omega^{n-1},\ 
\deg(p_\beta)\leq \frac{\deg(\omega)}{d}-A_\beta
\end{equation}) 
\end{prop}
This Proposition is  proved in \cite{mo} Proposition 1. 
The proof also gives us  the following algorithm  to find all the unknown 
data in the above equalities.

\begin{itemize}
\item
Input: The tame polynomial $f$ and $P\in\C[x]$ representing 
$[Pdx]\in H''$.

Output: $p_\beta,\beta\in I$ and $\xi$ satisfying (\ref{abelian})

We use the algorithm of \S \ref{homogeneous} and 
 write an element $\omega\in \Omega^{n+1}, \deg(\omega)=m$  in the form 
$$
\omega=\sum_{\beta\in I} p_\beta(g)\omega_\beta+dg\wedge 
d\psi,\ p_\beta\in\C[t],\ \psi\in \Omega^{n-1},\ 
\deg(p_\beta(g) \omega_\beta), \deg(dg\wedge d\psi) \leq m  
$$
This is possible because $g$ is homogeneous.
 We have
$$
\omega=\sum_{\beta\in I} p_\beta(f)\omega_\beta+df\wedge d\psi+\omega',\ 
\omega'=\sum_{\beta\in I} (p_\beta(g)-p_\beta(f))\omega_\beta+
d(g-f)\wedge d\psi
$$
The degree of $\omega'$ is strictly less than $m$ and so we repeat
what we have done at the beginning and finally we write $\omega$ as a 
$\C[t]$-linear combination of $\omega_\beta$'s. 
\end{itemize}
The algorithm for $H'$ is similar.
The statement about degrees is the direct consequence of the proof and
(\ref{14.11}). The procedure {\tt linear} takes a tame polynomial $f$  and 
$P\in\C[x]$  and returns a list. The first entry is a $1\times\mu$ matrix
$C=[p_\beta(t)]_{\beta\in I}$ and a $\mu\times \mu$ matrix representing a 
$(n-1)$-form $\xi$. 
The base ring must contain at  least one parameter and $C$ is written in 
the first parameter. The procedure {\tt linearp}  for $H'$ is similar to 
{\tt linear} for $H''$.
\section{Gauss-Manin connection}
\label{gaussmanin}
Let $S(t)\in\C[t]$ such that
$$
S(f)dx=df\wedge \eta_f,\ \eta_f=\sum_{i=1}^{n+1}(-1)^{i-1}p_i\widehat{dx_i}
\in \Omega^{n-1} 
$$
For instance one can take $S(t):=det(A_f-t.I)$, where $A_f$ is the 
multiplication by $f$ linear map form the Milnor vector space of $f$
$V:=\C[x]/\jacob(f)$ to itself. This definition of $S$ is implemented in the 
procedure ${\tt S}$ of {\tt brho.lib}.
 The Gauss-Manin connection associated to 
the fibration $f$ on $H''$ turns out to be a map 
$$
\nabla: H''\rightarrow H_C'', \nabla([Pdx])=\frac{[(Q_P-P.S'(f))dx]}{S}, \ 
P\in\C[x]
$$
where
\begin{equation}
\label{mande0}
Q_P=\sum_{i=1}^{n+1}
(\frac{\partial P}{\partial x_i}p_i+P\frac{\partial p_i}{\partial x_i})
\end{equation}
satisfying the Leibniz rule, where for a set $\tilde C\subset\C$ by 
$H_{\tilde C}''$ we mean the localization of $H''$ on the
multiplicative subgroup of $H''$ generated by $t-c,\ c\in \tilde C$.
Using the Leibniz rule one can extend $\nabla$ to a function from 
$H''_{C}$ 
to itself and  
so the iteration $\nabla^k=\nabla\circ\nabla\cdots\nabla$ $k$ times, 
makes sense. It is given by 
\begin{equation}
\label{mande1}
\nabla^k=
\frac{\nabla_{k-1}\circ\nabla_{k-2}\circ\cdots\circ\nabla_0}{S(t)^{k}}
\end{equation}
where
$$
\nabla_k:H''\rightarrow H'',\  \nabla_k([Pdx])=[(Q_P-(k+1)S'(t)P)dx]
$$
To calculate $\nabla:H'\rightarrow H'_C$ we use the fact that 
$$
\nabla^k \omega=\frac{\nabla^{k-1} d\omega}{df},\ \omega\in H'
$$
where $d: H'\rightarrow H''$ is taking differential and is well-defined. 
See \S 3 of \cite{mo} for more details on  $\nabla$.
Usually the iteration of the Gauss-Manin connection produces
polynomial forms with huge number of monomials. But fortunately our
Brieskorn module $H''$ (resp. $H'$) has already the canonical basis 
$\omega_\beta,\ \beta\in I$ (resp. $\eta_\beta,\ \beta\in I$) and after
writing $\nabla$ the obtained coefficients are much more easier to read.
The procedure ${\tt nabla}$ of {\tt brho.lib} uses the 
formulas (\ref{mande0}) and  (\ref{mande1}) and 
computes $\nabla$ and its iterations.  In $H''$ one can write
$$
S(t)\nabla(\omega_\beta)=\sum_{\beta'\in I}p_{\beta,\beta'}\omega_{\beta'},\ 
p_{\beta,\beta'}\in\C[t],\ \deg(p_{\beta,\beta'})\leq\deg(S)-1+
A_\beta-A_{\beta'}
$$
The bound on degrees can be obtained as follows: 
$$
S\omega_\beta=df\wedge \eta,\ \Rightarrow ds+dA_\beta=d+\deg(\eta)
$$
This is because $f$ is tame(see the proof of Lemma 4 of \cite{mo}).
$$
\deg(p_{\beta,\beta'})\leq \frac{\deg(d\eta)}{d}-A_{\beta'}
=s-1+A_\beta-A_{\beta'}
$$
The procedure {\tt nablamat} in {\tt brho.lib} calculates the matrix 
$\frac{1}{S(t)}[p_{\beta,\beta'}]$.
The Gauss-Manin connection $\nabla$ has two nice properties: 
\begin{enumerate}
\item
Griffiths transversality theorem: For all $i=1,2,\ldots,n+1$ we have 
$S(t)\nabla(F^i)\subset F^{i-1}$.
\item 
Residue killer: For all $\omega\in H$ there exists a $k\in \N$ such that
$\nabla^k\omega\in W_n$
\end{enumerate}
For the first one see \cite{gri83I}. 
The second one for $n=1$ is proved in
Lemma 2.3 of \cite{mov}. The proof for $n>1$ is similar and uses the
fact that the residue as a function in $t$ on a cycle around infinity is a
polynomial in $t$. 
\section{The numbers $d_\beta, \ \beta\in I$}
\label{dbeta}
Let $f$ be a tame polynomial with the last homogeneous part $g$, $F$ be its
homogenization and
$$
V=\C[x,x_0]/<\frac{\partial F}{\partial x_i}\mid i=1,2,\ldots,n+1>
$$ 
We consider $V$ as a $\C[x_0]$-module and it is shown in \cite{mo} that
$V$ is freely generated by $x^I:=\{x^\beta,\beta\in I\}$.
Let
$$
A_F:V\rightarrow V,\
A_F(G)=\frac{\partial F}{\partial x_0}G,\ G\in V
$$
\begin{prop}
The matrix of $A_F$ in the basis $x^I$ is
of the form $d.[x_0^{K_{\beta,\beta'}}c_{\beta,\beta'}]$,
where $K_{\beta,\beta'}:= d-1+\deg(x^\beta)-\deg(x^{\beta'})$ and $A_f:=[c_{\beta,\beta'}]$ 
is the multiplication by $f$ in the Milnor
vector space of $f$. In particular, if  $A_{\beta'}-A_{\beta}\geq 1$ then
$c_{\beta,\beta'}=0$ and 
$$\det (A_F-t.x_0^{d-1}I)=\det(A_f-t.I)
x_0^{(d-1)\mu}$$

\end{prop}
Using the above Proposition, the procedure {\tt muldF} 
calculates $A_F$.
\begin{proof}
Since the polynomial $F$ is weighted homogeneous, we have 
$\sum_{i=0}^{n+1}  \alpha_ix_i\frac{\partial F}{\partial x_i}=d.F$ and so 
$x_0\frac{\partial F}{\partial x_0}=d.F$ in $V$ (Note that $\alpha_0=1$ 
by definition). Let
\begin{equation}
\label{filmbardar}
F.x^\beta=\sum_{\beta'\in I}x^{\beta'}c_{\beta, \beta'}(x_0)+\sum_{i=1}^{n+1}
\frac{\partial F}{\partial x_i} q_i,\ c_{\beta,\beta'}(x_0)\in\C[x_0],\ 
q_i\in\C[x_0,x]
\end{equation}
Since the left hand side is homogeneous of degree $d+\deg(x^\beta)$ we 
can assume that the pieces of the write hand side are also homogeneous of
the same degree. This can be done by taking an arbitrary equation 
(\ref{filmbardar}) and subtracting the unnecessary parts. 
\end{proof}
Let $\tilde C$ be a finite subset of $\C$ and 
$\C[t]_{\tilde C}$ be the localization of
$\C[t]$ on its multiplicative subgroup generated by $t-c,\ c\in \tilde C$ 
and $F_t=F-t.x_0^{d}$. From now on we work with $\C[t]_{\tilde C}[x_0,x]$ 
instead of $\C[x_0,x]$ and redefine $V$ using  $\C[t]_{\tilde C}[x_0,x]$. 
Let
$$
V_{\tilde C}=\C[t]_{\tilde C}[x_0,x]/
<\frac{\partial F_t}{\partial x_0},\frac{\partial F}{\partial  x_i}, 
\mid i=1,2,\ldots,n+1>
$$
It is useful to reformulate $V_{\tilde C}$ in the following way: 
Let $R:=\C[t]_{\tilde C}[x_0]$ be the set of polynomials in $x_0$ with 
coefficients in $\C[t]_{\tilde C}$ and $A_t=A_F-t.d.x_0^{d-1}I$. We have
$$
V_{\tilde C}= V/ <\frac{\partial F_t}{\partial x_0}q\mid q\in V>=
R^\mu/A_t.R^\mu
$$
Here $R^\mu$ is the set of $\mu\times 1$ matrices with entries in $R$.
We consider  the statement:

$*(\tilde C)$: There is a function $\beta\in I\rightarrow
d_\beta\in \N\cup\{0\}$ such that the $\C[t]_{\tilde C}$-module $V_{\tilde C}$
is freely generated by 
\begin{equation}
\label{12.7.04}
\{x_0^{\beta_0}x^\beta, 0\leq \beta_0\leq d_\beta-1, \beta\in I\}
\end{equation}
To prove the statement $*(\tilde C)$ we may introduce a kind of Gaussian 
elimination in $A_t$ and simplify it. For this reason we introduce the
operation $GE(\beta_1,\beta_2,\beta_3)$.
For $\beta\in I$ let $(A_t)_\beta$ be the $\beta$-th row of $A_t$.
\begin{itemize}
\item
Input: $A_t$, $\beta_1,\beta_2,\beta_3\in I$ 
with $A_{\beta_1}\leq A_{\beta_2}$. 

Output: a matrix $A_t'$ and a finite subset $B$ of $\C$.

We replace $(A_t)_{\beta_2}$ with
$$
-\frac{(A_t)_{\beta_2,\beta_3}}{(A_t)_{\beta_1,\beta_2}}*(A_t)_{\beta_1} +
(A_t)_{\beta_2}
$$
and we set $B=zero(c(t))$, where 
$(A_t)_{\beta_1,\beta_2}=c(t).x_0^{K_{\beta_1,\beta_2}}$. 
Since  for all $\beta_4\in I$ we have
$$
K_{\beta_2,\beta_3}+K_{\beta_1,\beta_4}=K_{\beta_1,\beta_3}+
K_{\beta_2,\beta_4}
$$
The obtained matrix $A'_t$ is of the form 
$[x_0^{K_{\beta,\beta'}}c'_{\beta,\beta'}]$ and 
$c'_{\beta_2,\beta_3}=0$.  
If the matrix $B_t$ is obtained from $A_t$ by applying the above 
operation  and  $B\subset \tilde C$ then $A_t.R^\mu= B_tR^\mu$.
\end{itemize} 

We give  an example of algorithm which calculates $d_\beta$'s for 
for some finite set $\tilde C\subset \C$:
\begin{itemize}
\item 
Input: $A_t$ 

Output: $d_\beta,\beta\in I$ and a finite set $\tilde C\subset \C$

We identify $I$ with $\{1,2,\ldots,\mu\}$ and assume that 
$$
\beta_1\leq \beta_2 \Rightarrow A_{\beta_1}\geq A_{\beta_2} 
$$
The algorithm has $\mu$ steps indexed by $\beta=\mu,\mu-1,
\ldots,1$. We define the set $\tilde C$ to be empty.
In $\beta=\mu$ we have $A(\beta)=A_t$. In the step $\beta$  
we find the first $\beta_1$ such that $A(\beta)_{\beta,\beta_1}\not=0$ and 
put $d_{\beta_1}=d-1+\deg(x^\beta)-\deg(x^{\beta_1})$. 
For $\beta_2=\beta-1,\ldots,1$ we  make $GE(\beta,\beta_2,\beta_1)$ and define
$\tilde C=\tilde C\cup \cup_{\beta_2=1}^{\beta-1}B_{\beta_2}$, where 
$B_{\beta_2}$ is obtained  during $GE(\beta,\beta_2,\beta_1)$.
The numbers $d_\beta$'s obtained in this
way proves the statement $*(\tilde C)$.
\end{itemize}
The advantage of this algorithm 
is that in many cases it gives $\tilde C=C$. 
We do not have a proof
for $*(C)$.
One can also fix 
a value $c\in\C\backslash C$ and apply the above algorithm for $A_c$. In this
case we do note care about $\tilde C$ during the algorithm. 
The obtained $d_\beta$'s make the 
statement $*(\tilde C)$ true for some $\tilde C\subset \C$ with 
$c\not\in \tilde C$.  This algorithm is  implemented in the procedure
{\tt dbeta} of {\tt brho.lib}.
We prove the following weak statements:
\begin{prop}
There is a function $\beta\in I\rightarrow
d_\beta\in \N\cup\{0\}$ such that the $\C[t]_C$-module $V'$
is generated  by $\{
x_0^{\beta_0}x^\beta, 0\leq \beta_0\leq d_\beta-1, \beta\in I\}$.
\end{prop}
\begin{proof}
We have
$$
V'=R^\mu/A_t R^\mu
\stackrel{b}{\cong}
A_t^{-1} R^\mu/R^\mu= 
\frac{A_t^{{\mathrm adj}}R^\mu}{x_0^{\mu(d-1)}}/ R^\mu
$$
The isomorphism $b$ in the middle is obtained by acting $A_t^{-1}$
from left on $R^\mu$ and
${\mathrm adj}$ makes the adjoint of a matrix.
Now for $\beta\in I$ let $d_\beta$ be the pole order of $\beta$-th
arrow  of $\frac{A_t^{{\mathrm adj}}}{x_0^{\mu(d-1)}}$. The numbers $d_\beta$
are the desired numbers. It is easy to see that
$\{ x_0^{\beta_0}x^\beta, 0\leq \beta_0\leq d_\beta, \beta\in I\}$
generates $V'$.
\end{proof}
\begin{prop}
There is a subset $\tilde C\subset \C$ such that the statement 
$*(\tilde C)$ is true with $d_\beta=d-1,\beta\in I$.
\end{prop}
\begin{proof}
We identify $I$ with $\{1,2,\ldots,\mu\}$ and assume that 
$$
\beta_1\leq \beta_2 \Rightarrow A_{\beta_1}\geq A_{\beta_2} 
$$
By various use of operation $GE$ on $A_t$ we make all the entries 
of $(A_t)_{\beta,\mu}=0,\ \beta\in I\backslash \{\mu\}$. We repeat this
for $(A_t)_{\beta,\mu-1}=0, \beta\in I\backslash \{\mu,\mu-1\}$ and after
$\mu$-times we get a lower triangular matrix. We always divide on a polynomial
on $t$ with leading coefficient one and so  division by zero does not
occur. 
\end{proof}
\begin{prop}
Let $*(\tilde C)$ is valid with $d_\beta,\beta\in I$. Then
$$
A_\beta<n+1,\ d_\beta< d(n+2-A_\beta),\ \sum_{\beta\in I}d_\beta=\mu (d-1)
$$
\end{prop}
\begin{proof}
The first one is already in Steenbrink's Theorem \ref{stee}. The second 
inequality is obtained by applying the first inequality associated to 
$F-cx_0^d$ for some $c\in\C\backslash \tilde C$:
$$
A_{(d_\beta-1,\beta)}=A_\beta+\frac{d_\beta-1+1}{d}<n+2
$$ 
The Milnor number of $F-cx_0^d$ is $\sum_{\beta\in I}d_\beta$ and equals to  
the Milnor number of $g-cx_0^d$ which is $\mu (d-1)$. 
\end{proof}

\section{Main theorem of  \cite{mo}}
\label{maintheo}
Suppose that  $*(\tilde C)$ is valid with $d_\beta,\ \beta\in I$. 
Define
$$
I^k_{n+1}=\{\beta\in I\mid A_\beta= n+1-k\},\  
I^k_n=\{ \beta \in I \mid  A_\beta+\frac{1}{d}\leq n+1-k
\leq A_\beta +\frac{d_\beta}{d} \}
$$
\begin{theo}
For a tame polynomial $f$,  the set 
$$
B=\cup_{k=1}^{n}B^k_{n+1}\cup
 \cup_{k=0}^{n} B_{n}^k
$$
with
$$ 
B^{k}_{n+1}=\{\nabla^{n-k}\omega_\beta\mid  \beta\in I_{n+1}^k\},
B^{k}_{n}=\{\nabla^{n-k}\omega_\beta\mid  \beta\in I_n^{k}\},
$$
is a basis of $H''_{\tilde C}$ compatible with the mixed Hodge structure.
The same is true for $H'_{\tilde C}$ replacing 
$\nabla^{n-k} \omega_\beta$ with $\nabla^{n+1-k} \eta_\beta$. 
\end{theo}
Unfortunately, this theorem gives us a basis of a localization $H$ compatible
with mixed Hodge structure. In \S \ref{examples} we have computed such bases
for the Brieskorn module itself. 

To handle easier the pieces of the mixed Hodge structure of
$H_{\tilde C}$ we make the following table. 
\begin{center}
\begin{tabular}{|ccccccccccc|}
\hline
0 & & 1 & & 2 & & $\cdots$ & & $n$  &  & $n+1$   \\ 
\hline 
 & $I^{n}_n$ &$I_{n+1}^n$ & $I_{n}^{n-1}$ & $I_{n+1}^{n-1}$ & $I_n^{n-2}$  
& $\cdots$ & $I_n^{1}$ & 
$I_{n+1}^1$ &  $I_{n}^{0}$ &  \\
\hline
\end{tabular}
\end{center}
The procedure {\tt Imk} of {\tt brho.lib} gives us 
$x^\beta,\beta\in I_m^k,\ m=n,n+1, k=0,1,\ldots n$ with the order 
$I_n^n,I_n^{n-1},\ldots, I_n^0, I_{n+1}^n,I_{n+1}^{n-1},\ldots, I_{n+1}^1$. 
In the case $n=1$ we have the table
\begin{center}
\begin{tabular}{|ccccc|}
\hline
0 & & 1 & & 2   \\ 
\hline 
 & $I^1_1$ &$I_{2}^1$ & $I_1^0$ & \\
\hline
\end{tabular}
\end{center}
$$
I_1^1=\{\beta\in I\mid A_\beta+\frac{1}{d}\leq 1\leq A_\beta+
\frac{d_\beta}{d} \},\  
I_1^0=\{\beta\in I\mid A_\beta+\frac{1}{d}\leq 2\leq A_\beta+\frac{d_\beta}{d} \} 
$$
$$
I_2^1=\{\beta\in I\mid A_\beta=1 \}
$$
The forms $\omega_\beta,\ \beta\in I_1^1$ form a basis of $F^1\cap W_1$ and
the forms $\omega_\beta,\beta\in I_1^2$ form a basis of $H''/W_1$. 
Now to obtain a basis of $W_1/(F^1\cap W_1)$ we must modify $\nabla 
\omega_\beta,\ 
\beta\in I_1^0$.

The procedure {\tt changebase} calculates the matrix of the basis of the
Brieskorn module $H''_{\tilde C}$ obtained
in Theorem \ref{maintheo} in the canonical basis $\omega_\beta,\beta\in I$.
\section{Picard-Fuchs equations}
\label{pfeq}
It is a well-known fact that for a polynomial $f\in\C[x]$ and $\omega\in H$
the integral $I(t):=\int_{\delta_t} \omega$ satisfies
\begin{equation}
\label{shakira}
(\sum_{i=0}^kp_i(t)\frac{\partial^i}{\partial t^i})I_t=0,\ p_i(t)\in\C[t]
\end{equation}
called Picard-Fuchs equation, where $\delta_t\in H_n(L_t,\Z)$ is a
continuous family of topological cycles. 
When $f$ is tame, it is possible to calculate $p_i$' as follows:

We write 
$$
\nabla^i(\omega)=\sum_{\beta\in I}p_{i,\beta}\omega_\beta
$$
and define the $k\times \mu$ matrix $A=[p_{i,\beta}]$, where $i$ runs through
$1,2,\ldots,k$ and $\beta\in I$. Let $k$ be the smallest number such that
the the rows of $A_{k-1}$ are $\C(t)$-linear independent. Now, the rows of 
$A_{k}$ are $\C(t)$-linear dependent and this gives us (after multiplication
by a suitable element of $\C[t]$) 
$$
\sum_{i=0}^{k}p_i(t)\nabla^i(\omega)=0,\ p_i(t)\in\C[t]
$$
Using the formula (\ref{gofta}) and integrating the above equality,
we get the equation (\ref{shakira}).  
The procedure {\tt PFeq} from the library {\tt brho.lib} calculates
$p_i$'s in (\ref{shakira}). 
\section{Polynomials in one variable, $n=0$}
\label{n=0}
The theory developed in \S \ref{adami} does 
not work for the case $n=0$. 
For a polynomial of degree $d$ in one variable 
$\dim (H^0(L_t,\C))=d$ but $\mu=d-1$. However, if we use the following
definition of homology and cohomology for a discrete topological space 
$M$, $$
H_0(M,\Z)=\{m=\sum_{i}a_im_i\mid a_i\in\Z,\ m_i\in M\mid 
deg(m)=\sum_i a_i=0\}
$$
$$
H^0(M,\C)=\{f: H_0(M,\Z)\rightarrow \C \hbox{ linear}\}/\{f\mid \hbox{ f is constant on } M\}
$$
then
$$
H'=\C[x]/\C[f], \ H''=\C[x]dx/f'\C[f]dx, \ I=\{1,x,x^2,\ldots,x^{d-2}\}, 
\mu=d-1
$$
In this case 
$$
\int_\delta\omega =\sum_{i} a_i\omega(p_i), \
\hbox{ where } \delta=\sum_{i}a_ip_i,\ a_i\in\Z, \ p_i\in f^{-1}(t), \ 
\omega\in H'
$$
If, for instance, $f'=0$ has $d$ distinct root then every vanishing cycle
in  $L_t$ is a difference of two points of $L_t$. 
The set $B=\{x,x^2,\ldots, x^{d-1}\}$ form a basis of $H'$ and its 
$\nabla$ which is  $\{dx,xdx,\ldots x^{d-2}dx\}$ 
(up to multiplication by some 
constants) form a basis of $H''$. The first fact is easy to see. We write
$f= a_{d}x^{d}+f_0$ and for a polynomial 
$p(x)\in\C[x]$ whenever we find some
$x^{d}$ we replace it with $\frac{f-f_0}{a_{d}}$ and at the end we get 
$p(x)=p_0(f)+\sum_{i=1}^{d-1}p_i(f)x^i$ or equivalently 
$p=\sum_{i=1}^{d-1} p_i(t)x^i$ in $H'$. There is no $\C[t]$-linear 
relation between the elements of $B$ because $B$ restricted to each regular 
fiber is of dimension $d$.
We write 
$$
p(x)dx=\sum_{i=0}^{d-2 }q_i(f)x^idx+q_{d-1}(f)x^{d-1}dx=
(\sum_{i=0}^{d-2 }q_i(f)x^idx-\frac{q_{d-1}(f)f_0'}{d.a_d}dx)+
\frac{q_{d-1}(f)f'}{d.a_{d}}dx
$$
and this proves the statement for $H''$.

The proposition (\ref{4.7.03}) can be stated in the case $n=0$ as follows: 
The only case in which $dA_\beta-\alpha_i=0$ is when $n=0$ and $P=1$.
In the case $n=0$ for $P\not=1$ we have
$$
\frac{\partial g}{\partial x_i}.Pdx=\frac{d}{d.A_\beta-\alpha_i}
\frac{\partial P}{\partial x_i}gdx
$$ 
and if $P=1$ then $\frac{\partial g}{\partial x_i}.Pdx$ is zero in $H''$.
The argument in (\ref{29.6}) and (\ref{roesler}) can be done also in 
the case $n=0$. In this case if
\begin{equation}
\eta=\sum_{\beta\in I}p_\beta(g)\eta_\beta +p(g),\ p, p_\beta\in
\C[t]  
\end{equation}
where each piece in the right hand side of the above equality has degree less 
than $\deg(\eta)$  then
\begin{equation}
d\eta=\sum_{\beta\in I}(p_\beta(g)A_\beta+p'_\beta(g)g)\omega_\beta+
p'(g)dg
\end{equation}
Based on this observation, the procedure {\tt linear, linearp,} works
for the case $n=0$.   
 
In the case $n=0$, we have only the set  
$I_0^0=\{A_\beta+\frac{1}{d}\leq 1\leq A_\beta+\frac{d_\beta}{d}\}$ and this 
is equal to $I$.  
We have $d_\beta<d.(n+2-A_\beta)=2d-\beta-1=$ and $A_\beta=\frac{\beta+1}{d}$.
We conclude that 
$$
d\leq d_\beta+\beta+1<2d
$$
Now the infinitesimal Hilbert problem (see \cite{Il02} Problem 7) can be 
stated in the case $n=0$. Can one give an effective solution to this problem
in this case? The positive answer to this question may give light into the the 
problem in the case $n=1$.
\section{Examples}
\label{examples}
{\tiny 
\begin{multicols}{2}
For all the examples bellow we run
\begin{verbatim}
>LIB "brho.lib";
>LIB "matrix.lib";
\end{verbatim}
\subsection{Examples, $n=0$}
For examples of this section we run
\begin{verbatim}
>ring r0=(0,t),x, dp;
\end{verbatim}
\begin{exam}\rm
$f=x^5-5x$, $P=\{\epsilon^i\mid i=0,1,2,3\}$, 
$C=\{-4\epsilon^i\mid i=0,1,2,3\}$, where 
$\epsilon=e^{\frac{2\pi i}{d-1}}$ is the $d$-th root 
of unity.
\begin{verbatim}
> int d=5; poly f=x^d-d*x; okbase(std(jacob(f)));
_[1]=x3
_[2]=x2
_[3]=x
_[4]=1
> Abeta(f);
_[1,1]=4/5
_[2,1]=3/5
_[3,1]=2/5
_[4,1]=1/5 
> poly Sf=S(f); Sf;
(t4-256)
> list l=nablamat(f,Sf);
> l[1]; print(l[2]);
1/(5t4-1280)
(-t3), 128,   (-48t),(16t2),
(4t2), (-2t3),192,   (-64t),
(-16t),(8t2), (-3t3),256,   
64,    (-32t),(12t2),(-4t3)
//This is the matrix of nabla in the canonical 
//basis x^3,x^2,x^1,1.
 >PFeq(f,1);
_[1,1]=6144
_[1,2]=(35625t)
_[1,3]=(33375t2)
_[1,4]=(8750t3)
_[1,5]=(625t4-160000)
\end{verbatim}
The  residues of $\frac{dx}{f-t}$ at its poles satisfy the Picard-Fuchs
equation
$$
6144+35625t\frac{\partial}{\partial t}+33375t^2\frac{\partial^2}{\partial t^2}
8750t^3\frac{\partial^3}{\partial t^3}+ 
$$
$$
(625t^4-160000)
\frac{\partial^4}{\partial t^4}=0
$$
\end{exam}
\subsection{Examples $n=1$}
For the examples bellow we define
\begin{verbatim}
  ring r1=(0,t), (x,y), dp;
\end {verbatim}

\begin{exam}\rm 
$f=xy(x+y-1)$.
\begin{verbatim}
 > poly f= x2y+xy2-xy ;
 > poly g=lasthomo(f); g;
x2y+xy2
 > okbase(std(jacob(g)));
_[1]=y2
_[2]=y
_[3]=x
_[4]=1
 > print(muldF(f-par(1)));
(-3t+1/18)*x2,-1/18*x3,    0,       0,      
1/6*x,        (-3t-1/6)*x2,0,       0,      
1/6*x,        -1/6*x2,     (-3t)*x2,0,      
1/2,          -1/2*x,      0,       (-3t)*x2
 > poly Sf=S(f); Sf;
(t4+1/27t3)
//We can take Sf=t*(t+1/27);
 > list l1=nablamat(f,Sf);
 > l1[1]; " "; print(l1[2]);
1/(54t2+2t)
 
(18t+1),(-18t-1),0,(-2t),
1,      -1,      0,(-6t),
1,      -1,      0,(-6t),
3,      -3,      0,(-18t)
//--------------
 > dbeta(f,par(1));
0,2,2,4 
 > Imk(f,par(1));
[1]:
   [1]:
      [1]:
         1
   [2]:
      [1]:
         1
[2]:
   [1]:
      [1]:
         x
      [2]:
         y
 > list l3=changebase(f,Sf,par(1));
 > print(l3[1]); " "; print(l3[2]); det(l3[2]);
1,3/(54t2+2t),1,1
 
0,0, 0,1,    
1,-1,0,(-6t),
0,0, 1,0,    
0,1, 0,0     
1
//--------------
 > dbeta(f);
2,2,2,2 
 > Imk(f);
[1]:
   [1]:
      [1]:
         1
   [2]:
      [1]:
         y2
[2]:
   [1]:
      [1]:
         x
      [2]:
         y
 > list l2=changebase(f,Sf);
 > print(l2[1]); " "; print(l2[2]); det(l2[2]);
1,1/(54t2+2t),1,1
 
0,      0,       0,1,    
(18t+1),(-18t-1),0,(-2t),
0,      0,       1,0,    
0,      1,       0,0     
(18t+1)
//The obtained basis does not work for the fiber c=-1/18
//--------------
 > PFeq(f,1, Sf);
_[1,1]=6
_[1,2]=(54t+1)
_[1,3]=(27t2+t)
_[1,4]=0
_[1,5]=0
\end{verbatim}
We get the following basis of $H''$ compatible with mixed Hodge structure.
\begin{center}
\begin{tabular}{|c|c|}
\multicolumn{2}{c}{$f=xy(x+y-1)$}
\\ \hline
$\gr^1_F\gr_{1}^W H''$ & $[1]$ \\
\hline 
$\gr^0_F\gr_{1}^WH'' $ &  $[y^2]-[y]-6t[1]$ \\
 
\hline
$\gr^1_F\gr_{2}^WH ''$ &  $[x],[y]$ \\
\hline
\end{tabular}
\end{center}
The integrals $I=\int_{\delta_t}\frac{\omega}{f-t}$ satisfy the Picard-Fuchs equation
$$
6+(54t+1)\frac{\partial I}{\partial t}+(27t^2+t)\frac{\partial^2 I}{\partial t^2}=0
$$
\end{exam}
\begin{exam}\rm
$f=2(x^3+y^3)-3(x^2+y^2)$, 
$P=\{(0,0),(0,1),(1,0),(1,1)\}$, 
$C=\{0,-1, -1, -2\}$, 
\begin{verbatim}
 > poly f= 2*x3+2*y3-3*x2-3*y2 ;
 > poly g=lasthomo(f); g;
2*x3+2*y3
 > okbase(std(jacob(g)));
_[1]=xy
_[2]=y
_[3]=x
_[4]=1
 >S(f);
(t4+4t3+5t2+2t)
//We can put
 >poly Sf=t*(t+1)*(t+2);
 > list l2=changebase(f,Sf);
 > print(l2[1]); " "; print(l2[2]); det(l2[2]);
1,-1/(6t+12),1,1
 
0, 0,0,1,
-2,1,1,0,
0, 0,1,0,
0, 1,0,0 
-2
\end{verbatim}
\begin{center}
\begin{tabular}{|c|c|}
\multicolumn{2}{c}{$f=2(x^3+y^3)-3(x^2+y^2)$}
\\ \hline
$\gr^1_F\gr_{1}^W H''$ & $[1]$ \\
\hline 
$\gr^0_F\gr_{1}^WH'' $ &  $[2xy-x-y]$ \\
 
\hline
$\gr^1_F\gr_{2}^WH'' $ &  $[x],[y]$ \\
\hline
\end{tabular}
\end{center}
\end{exam}
\begin{exam}\rm
$f=x^4+y^4-x$.
\begin{verbatim}
 > poly f= x4+y4-x ;
 > poly g=lasthomo(f);
 > okbase(std(jacob(g)));
_[1]=x2y2
_[2]=xy2
_[3]=x2y
_[4]=y2
_[5]=xy
_[6]=x2
_[7]=y
_[8]=x
_[9]=1
 > poly Sf=S(f); Sf;
(t9+81/256t6+2187/65536t3+19683/16777216)
//We can take 
 >Sf=t^3+27/256;
 > dbeta(f,par(1));
2,2,2,5,2,2,5,2,5 
 > Imk(f,par(1));
[1]:
   [1]:
      [1]:
         1
      [2]:
         x
      [3]:
         y
   [2]:
      [1]:
         y
      [2]:
         y2
      [3]:
         x2y2
[2]:
   [1]:
      [1]:
         x2
      [2]:
         xy
      [3]:
         y2
 > list l3=changebase(f,Sf,par(1));
 > print(l3[1]); " "; print(l3[2]); det(l3[2]); 
1,1,1,4/(256t3+27),24/(256t3+27),1/(256t3+27),1,1,1
 
0,      0,    0,0,0,0,0,      0,1,
0,      0,    0,0,0,0,0,      1,0,
0,      0,    0,0,0,0,1,      0,0,
0,      0,    9,0,0,0,(-16t2),0,0,
3,      (-2t),0,0,0,0,0,      0,0,
(128t2),9,    0,0,0,0,0,      0,0,
0,      0,    0,0,0,1,0,      0,0,
0,      0,    0,0,1,0,0,      0,0,
0,      0,    0,1,0,0,0,      0,0 
(2304t3+243)  // 9*256*Sf;
 > matrix A=l3[2];
 > A[6,1..ncols(A)]= 
((-128*t2)/3)*submat(A,5,1..ncols(A))+submat(A,6,1..ncols(A));
 > A[5,1..ncols(A)]= 
2*t*submat(A,6,1..ncols(A))+submat(A,5,1..ncols(A)); print(A);
 0,0,0,0,0,0,0,      0,1,
 0,0,0,0,0,0,0,      1,0,
 0,0,0,0,0,0,1,      0,0,
 0,0,9,0,0,0,(-16t2),0,0,
 1,0,0,0,0,0,0,      0,0,
 0,1,0,0,0,0,0,      0,0,
 0,0,0,0,0,1,0,      0,0,
 0,0,0,0,1,0,0,      0,0,
 0,0,0,1,0,0,0,      0,0
\end{verbatim}
We obtain the following table
\begin{center}
\begin{tabular}{|c|c|}
\multicolumn{2}{c}{$f=x^4+y^4-x$}
\\ \hline
$\gr^1_F\gr_{1}^W H''$ & $[1],[x],[y]$ \\
\hline 
$\gr^0_F\gr_{1}^WH'' $ &  $9[x^2y]-16t^2[y], [x^2y^2], [xy^2]$ \\
 
\hline
$\gr^1_F\gr_{2}^WH'' $ &  $[x^2], [xy],[y^2]$ \\
\hline
\end{tabular}
\end{center}
We make the following remark
\begin{verbatim}
 > reduce(9*x2*y-16*(f^2)*y, std(jacob(f)));
0
\end{verbatim}
\end{exam}
\end{multicols} 
} 

\def\cprime{$'$} \def\cprime{$'$}

\bibliographystyle{plain}

\end{document}